\documentclass[11pt,a4paper]{article}

\usepackage{ifpdf}

\makeatletter
\let\@font@warningori\@font@warning
\newcommand\shutup{\def\@font@warning##1{}}
\newcommand\youcanspeak{\let\@font@warning\@font@warningori}
\makeatother

\usepackage{amsmath,amsbsy,bm,latexsym,enumerate}

\shutup
\usepackage{fourier}
\youcanspeak
\usepackage[scaled=0.875]{helvet}

\usepackage[utf8x]{inputenc}
\usepackage[german,dutch,english]{babel}

\usepackage{amsthm}

\usepackage{graphicx,epsfig,color}

\usepackage{tikz}
\usetikzlibrary{calc}
\usetikzlibrary{intersections}

\usepackage{anysize}

\usepackage{subfigure}

\newcommand {\cA}         {\mathcal{A}}
\newcommand {\cB}         {\mathcal{B}}
\newcommand {\cV}         {\mathcal{V}}

\renewcommand{\epsilon}{\varepsilon}

\renewcommand{\leq}{\leqslant}
\renewcommand{\geq}{\geqslant}

\theoremstyle{plain}
\newtheorem{theorem}{Theorem}
\newtheorem{lemma}[theorem]{Lemma}

\theoremstyle{definition}
\newtheorem{definition}[theorem]{Definition}

\title{Constructing subset partition graphs with strong adjacency and end-point count properties}
\author{Nicolai Hähnle\\\texttt{haehnle@math.tu-berlin.de}}
\date{\today}

\begin{document}

\maketitle

\begin{abstract}
  Kim defined a very general combinatorial abstraction
  of the diameter of polytopes called subset partition graphs
  to study how certain combinatorial properties of such graphs may be achieved in lower bound constructions.
  Using Lovász' Local Lemma,
  we give a general randomized construction for subset partition graphs
  satisfying strong adjacency and end-point count properties.
  This can be used as a building block to conceptually simplify the constructions
  given in~\cite{arXiv:1103.3362v1}.

  We also use our method to construct abstract spindles,
  an analogy to the spindles used by Santos~\cite{Santos10} to disprove
  the Hirsch conjecture,
  of exponential length
  which satisfy the adjacency and end-point count properties.
\end{abstract}

\section{Introduction}

One of the great mysteries of convex geometry is the question
whether the diameter of the vertex-edge graph of polyhedra is always bounded
by a polynomial in the dimension and the number of facets.
A positive answer to this question is necessary for the existence of a polynomial-time pivot rule for
the Simplex method.
Yet the best known upper bound is $n^{1 + \log d}$,
where $n$ is the number of facets and $d$ is the dimension~\cite{MR1130448}.
In the other direction, the best known constructions provide a diameter
of $(1 + \varepsilon)n$ for some fixed $\varepsilon > 0$
in sufficiently high but constant dimension~\cite{Santos10},
so a considerable gap in our knowledge remains.
For a recent survey of the field, see~\cite{MR2681516}.

There is a long history of
research considering purely combinatorial abstractions of the graphs of polytopes.
Important representatives of this approach can be found in~\cite{MR0496673, MR1176376}.
This has ultimately lead to a combinatorial abstraction in which
the best known general upper bounds can be proved,
but which also admits constructions of at least an almost quadratic diameter~\cite{EHRR10}.

Kim~\cite{arXiv:1103.3362v1} suggests a more systematic study
of the various combinatorial properties of graphs of polytopes.
The goals of this endeavor are to identify properties that would help establish a polynomial
upper bound on the diameter,
and to investigate combinatorial constructions with high diameter,
possibly as an inspiration for more geometric constructions.
To this end, he gives the following definition:
\begin{definition}
  Let $S$ be a set of \emph{symbols} and $n = |S|$.
  Let $\cA \subseteq \binom{S}{d}$, where $1 \leq d \leq n$ is called the \emph{dimension}.
  Let $G = (\cV, E)$ be a connected undirected graph where the vertices $V \in \cV$
  form a partition of $\cA$.
  In particular, every $V \in \cV$ is a non-empty subset of $\cA$.
  Then we call $G$ a \emph{$d$-dimensional subset partition graph of $\cA$ on the symbol set $S$}.
\end{definition}
The symbols correspond to the facets of a polyhedron and $d$ to its affine dimension,
so that the sets $A \in \cA$ correspond to non-degenerate vertices,
which can be identified by the $d$ facets in which they are contained.
Note that we can restrict ourselves to simple polyhedra,
i.e. polyhedra without degenerate vertices,
by a perturbation argument, see~\cite{Santos10} for a modern treatment.

By itself, the definition of subset partition graphs provides very little structure.
However, it can easily be augmented with additional properties.
In the search for such properties,
we are guided by what is known to hold for the graphs of polyhedra.
Mostly following~\cite{arXiv:1103.3362v1}, we say that $G$ satisfies
\begin{itemize}
  \item \emph{adjacency}, if $A, A' \in \cA$ and $|A \cap A'| = d-1$ implies that $A$ and $A'$
    are in the same or adjacent vertices of $G$;
  \item \emph{strong adjacency}, if, in addition to adjacency, for every pair of adjacent vertices $V, V' \in \cV$
    one has $A \in V$ and $A' \in V'$ such that $|A \cap A'| = d-1$;
  \item \emph{end-point count}, if for all $F \in \binom{S}{d-1}$
    one has $\#\{ A \in \cA : A \supset F \} \leq 2$; and
  \item the \emph{one-subset} or \emph{singleton property}, if $|V| = 1$ for all $V \in \cV$.
\end{itemize}
Furthermore, we define the \emph{restriction} of $G$ with respect to a subset $F \subset S$
as the subset partition graph $G|_F = (\cV|_F, E|_F)$, where
\begin{align*}
 \cV|_F &= \{ V|_F ~:~ V \in \cV \text{ and } V|_F \text{ not empty} \}, \text{ where} \\
   V|_F &= \{ A \setminus F ~:~ A \in V \text{ and } A \supseteq F \} \text{ for all } V \in \cV,
\end{align*}
and the adjacency structure $E|_F$ is the obvious induced graph.
Then we say that $G$ satisfies
\begin{itemize}
  \item \emph{dimension reduction}, if $G|_F$ is connected (or empty) for every $F \in \binom{S}{\leq d}$.
\end{itemize}
Kim points out that this provides a framework in which previous abstractions can be understood
by requiring that some subset of these properties is satisfied.
For example, the connected layer families of~\cite{EHRR10}
are subset partition graphs where the graph structure is a path,
and dimension reduction is satisfied.
For the graph of every simple polyhedron,
there is a natural corresponding subset partition graph
that satisfies all these properties.

Motivated by Santos' counter-example to the Hirsch conjecture~\cite{Santos10},
he also defines what we may call \emph{abstract spindles} or \emph{abstract Dantzig figures}.
Suppose $n = 2d$ and we distinguish two \emph{apices} $A_1, A_2 \in \cA$ such that $A_1 \cup A_2 = S$,
then we may call $G$ an abstract spindle.
Its \emph{length} is the distance between the vertices $V_1, V_2 \in \cV$
that contain $A_1$ and $A_2$, respectively.

The contribution of this paper is two-fold.
First, we describe a randomized construction, using Lovász' Local Lemma,
    that turns any SPG that satisfies the singleton property
    into an SPG that additionally satisfies strong adjacency and end-point count.
    This transformation is at the cost of multiplying the dimension and number of symbols by a factor
    that depends on the maximum degree of the underlying graph.
    It can be used as a building block for constructions.
    In particular, it can be used to conceptually simplify constructions of~\cite{arXiv:1103.3362v1},
    as long as one does not care about constants.
We then construct spindles of exponential length satisfying
    both strong adjacency and end-point count properties.
    Previously, only a polynomial length construction was given in~\cite{arXiv:1103.3362v1}.

\section{Obtaining the strong adjacency and end-point count properties}

Let $G = (\cV, E)$ be a $d$-dimensional subset partition graph on the symbol set $S$
that satisfies the singleton property.
Our goal is to use this graph as a template for the construction of an
$rd$-dimensional SPG $G' = (\cV', E')$ on the symbol set $S' = S \times [r]$ such that:
\begin{enumerate}
  \item As an undirected graph, $G'$ is isomorphic to a subdivision of $G$.
  \item $G'$ satisfies the singleton property.
  \item Let $V \in \cV$, $V = \{ A \}$, and let $V' \in \cV'$ the corresponding
    vertex under the mentioned isomorphism. Then $V' = \{ A \times [r] \}$.
  \item $G'$ satisfies the strong adjacency and end-point count properties.
\end{enumerate}
The multiplier $r$ is a parameter that we would like to choose as small as possible.
We will give a randomized construction that shows how $r$ can be bounded
in terms of the maximum degree $\Delta$ of $G$.
Our construction has two steps.
In the first step, we replace every set $A \in \cA$ by $A \times [r]$
to satisfy the adjacency and end-point count properties.
Clearly, we must choose $r \geq 2$ for this to work.
In the second step, we subdivide edges of $G$ to ``interpolate''
between the sets in the end-point vertices and thus establish strong adjacency.
Locally, the construction will look like the following example:
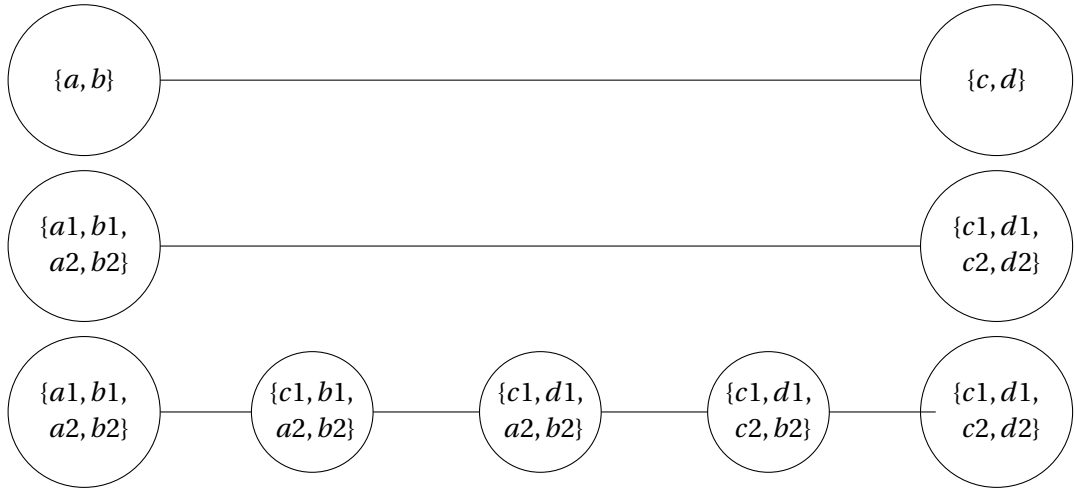
\begin{center}
  \begin{tikzpicture}
    \begin{scope}
    \draw (0,0) circle (1cm) node {$\{a,b\}$};
    \draw (12,0) circle (1cm) node {$\{c,d\}$};
    \draw (1,0) -- (11,0);
    \end{scope}

    \begin{scope}[yshift=-2.2cm]
    \draw (0,0) circle (1cm) node {$\begin{array}{r@{\,}l}\{a1,&b1,\\a2,&b2\}\end{array}$};
    \draw (12,0) circle (1cm) node {$\begin{array}{r@{\,}l}\{c1,&d1,\\c2,&d2\}\end{array}$};
    \draw (1,0) -- (11,0);
    \end{scope}

    \begin{scope}[yshift=-4.4cm]
    \draw (0,0) circle (1cm) node {$\begin{array}{r@{\,}l}\{a1,&b1,\\a2,&b2\}\end{array}$};
    \draw (3,0) circle (0.8cm) node {$\begin{array}{r@{\,}l}\{c1,&b1,\\a2,&b2\}\end{array}$};
    \draw (6,0) circle (0.8cm) node {$\begin{array}{r@{\,}l}\{c1,&d1,\\a2,&b2\}\end{array}$};
    \draw (9,0) circle (0.8cm) node {$\begin{array}{r@{\,}l}\{c1,&d1,\\c2,&b2\}\end{array}$};
    \draw (12,0) circle (1cm) node {$\begin{array}{r@{\,}l}\{c1,&d1,\\c2,&d2\}\end{array}$};
    \draw (1,0) -- (2.2,0) (3.8,0) -- (5.2,0) (6.8,0) -- (8.2,0) (9.8,0) -- (11.2,0);
    \end{scope}
  \end{tikzpicture}
\end{center}
The following alternative illustration of the same situation,
inspired by~\cite{arXiv:1103.3362v1},
represents the sets of $\cA$ as subsets of squares of a fixed-size grid.
The symbols of $S$ are represented as a row of squares.
In the first step of the construction, this row is replicated $r$ times (here, $r = 2$).
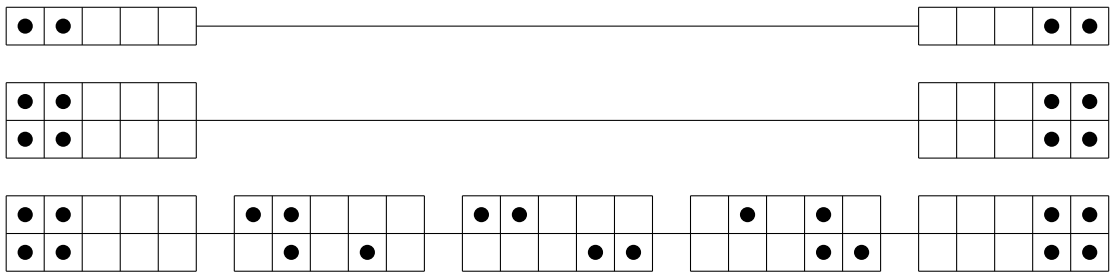
\begin{center}
  \def\drawgrid#1#2#3#4{
    \begin{scope}[shift={#3}]
    \begin{scope}[scale=0.5]
      \foreach \x in {0,1,...,#1}
        \draw (\x,0) -- (\x,#2);
      \foreach \y in {0,1,...,#2}
        \draw (0,\y) -- (#1,\y);
      \foreach \ofs in {#4}
        \fill \ofs +(0.5,0.5) circle (0.2);
    \end{scope}
    \end{scope}
  }

  \begin{tikzpicture}
    \begin{scope}
      \drawgrid{5}{1}{(0,0)}{(0,0),(1,0)};
      \drawgrid{5}{1}{(12,0)}{(3,0),(4,0)};
      \draw (2.5,0.25) -- (12,0.25);
    \end{scope}

    \begin{scope}[shift={(0,-1.5)}]
      \drawgrid{5}{2}{(0,0)}{(0,0),(1,0),(0,1),(1,1)};
      \drawgrid{5}{2}{(12,0)}{(3,0),(4,0),(3,1),(4,1)};
      \draw (2.5,0.5) -- (12,0.5);
    \end{scope}

    \begin{scope}[shift={(0,-3)}]
      \drawgrid{5}{2}{(0,0)}{(0,0),(1,0),(0,1),(1,1)};
      \drawgrid{5}{2}{(3,0)}{(3,0),(1,0),(0,1),(1,1)};
      \drawgrid{5}{2}{(6,0)}{(3,0),(4,0),(0,1),(1,1)};
      \drawgrid{5}{2}{(9,0)}{(3,0),(4,0),(3,1),(1,1)};
      \drawgrid{5}{2}{(12,0)}{(3,0),(4,0),(3,1),(4,1)};

      \draw (2.5,0.5) -- (3,0.5) (5.5,0.5) -- (6,0.5) (8.5,0.5) -- (9,0.5) (11.5,0.5) -- (12,0.5);
    \end{scope}
  \end{tikzpicture}
\end{center}
In the example, we interpolated between the endpoints one row at a time,
so that for every intermediate vertex, at most one row differs from the corresponding row
at both endpoints.
This property will be useful for our construction as well.

So here is the second step of our construction.
For every edge between a pair of vertices $V = \{ A \}$, $V' = \{ A' \}$,
we choose a permutation $\pi$ of the $r$ rows uniformly at random.
We subdivide the edge $VV'$ into a path of length $r \cdot (d - |A \cap A'|)$.
We consider this path as composed of $r$ segments of length $d - |A \cap A'|$ each.
On every edge of the $j$-th segment, we take an element of $A \setminus A'$ from the $\pi(j)$-th row
and replace it with an element of $A' \setminus A$ from the same row.
The order in which elements within a single row are replaced is arbitrary;
what matters is that rows are adjusted one at a time,
and that the order of the rows is chosen uniformly at random,
independently for each edge.

Since we feel that it would not help understanding,
we do not provide a complete formalization of this construction.
We trust that the reader may develop one herself if she deems it necessary.
We only mention explicitly that the $\pi(j)$-th row is the set $S \times \{ \pi(j) \} \subset S'$.
For a set $B \subseteq S' = S \times [r]$ we use the notation
\[ B|_{\pi(j)} := \{ x \in S ~:~ (x,\pi(j)) \in B \} \]
for the $\pi(j)$-th row of $B$.

It is easy to see that the resulting subset partition graph $G'$
satisfies the first three properties demanded at the beginning of this section.
Furthermore,
for every pair of adjacent vertices $V = \{ A \}$, $V' = \{ A' \}$ one has
$|A \cap A'| = rd - 1$ by construction,
and so the second part of the definition of strong adjacency is also satisfied.
It remains to show that $G'$ satisfies adjacency and the end-point count property with positive probability.
\begin{lemma}
  \label{lemma:adjacency}
  If $r \geq \lceil 16e\Delta \rceil$,
  the following holds with positive probability:
  Let $A, A' \in \cA'$ such that $|A \cap A'| = rd - 1$.
  Then $A$ and $A'$ are contained in adjacent vertices of $G'$.
\end{lemma}
This Lemma, once proved, clearly establishes adjacency.
It also implies the end-point count property,
because $G'$ satisfies the singleton property and
contains no triangles (every edge of $G$ is subdivided into a path of length at least $r$).
The remainder of this section is devoted to a proof of Lemma~\ref{lemma:adjacency}
using Lovász' Local Lemma.
We begin with a restatement of an earlier observation:
\begin{lemma}
  \label{lemma:subdivided-vertex-rows}
  Let $A \in \cA'$ and let $e \in E$ be an edge of $G$ on whose subdivision $A$ lies
  ($e$ is only guaranteed to be unique if $A$ does not correspond to a vertex of $G$).
  Let $V, V' \in \cV$ be the end-points of $e$.
  Then all but at most one row of $A$ are equal to either $V$ or $V'$.
\end{lemma}

\begin{lemma}[Localization Lemma]
  Suppose $r \geq 4$.
  Let $A_1, A_2 \in \cA'$ such that $|A_1 \cap A_2| = rd - 1$
  and let $e_1, e_2 \in E$ be edges of $G$ on whose subdivision $A_1$ and $A_2$ lie, respectively.
  Then $e_1$ and $e_2$ share an end-point.
\end{lemma}
\begin{proof}
  Our goal is to show that there is one row in which $A_1$ and $A_2$ coincide,
  and which also coincides with both an end-point of $e_1$ and an end-point of $e_2$.
  Then these end-points must in fact be the same vertex, that is,
  $e_1$ and $e_2$ share an end-point.

  There is one row which must be ruled out because $A_1$ and $A_2$ differ in exactly one row.
  In addition, up to one row must be ruled out since, by Lemma~\ref{lemma:subdivided-vertex-rows},
  $A_1$ might not coincide with either end-point of $e_1$ in one of the rows.
  Similarly, one additional row might have to be ruled out due to $A_2$.
  In total, up to three rows must be eliminated from consideration.
  Since $r \geq 4$, there remains at least one row in which $A_1$ and $A_2$ coincide,
  and which also coincides with one end-point of both edges $e_1$ and $e_2$.
  This completes the proof.
\end{proof}

The Localization Lemma suggests the following definition.
Let $V \in \cV$ and let $e_1, e_2 \in E$ edges incident to $V$.
Let furthermore $V' \in \cV'$ be the vertex in $G'$ corresponding to $V$.
Then the \emph{bad event} $B_{e_1 e_2}$ is the event that
there exist vertices $W_1 = \{ B_1 \}, W_2 = \{ B_2 \} \in \cV' \setminus \{ V' \}$
on the subdivisions of $e_1$ and $e_2$, respectively,
with $|B_1 \cap B_2| = rd - 1$.
If none of these bad events occur, then by the Localization Lemma,
the constructed $G'$ has all desired properties.

\begin{figure}
  \begin{center}
  \begin{tikzpicture}
    \begin{scope}
      \node[draw,shape=circle] (V) at (0,0) {$V\phantom{'}$};
      \node[draw,shape=circle] (V1) at (4,1.5) {$V_1\phantom{'}$};
      \node[draw,shape=circle] (V2) at (4,-1.5) {$V_2\phantom{'}$};

      \draw (V) -- node[above] {$e_1$} (V1);
      \draw (V) -- node[below] {$e_2$} (V2);
    \end{scope}

    \begin{scope}[shift={(6,0)}]
      \node[draw,shape=circle] (Vp) at (0,0) {$V'$};
      \node[draw,shape=circle] (Vp1) at (4,1.5) {$V_1'$};
      \node[draw,shape=circle] (Vp2) at (4,-1.5) {$V_2'$};

      \draw (Vp) -- node[above] {$e_1$} (Vp1);
      \draw (Vp) -- node[below] {$e_2$} (Vp2);

      \foreach \x in {0.8,1.2,...,3.2}
        \fill ($(\x,0) + {\x}*(0,0.375)$) circle (2pt);

      \fill ($(2.4,0) + 2.4*(0,0.375)$) circle (3pt) node [below right] {$W_1$};

      \foreach \x in {0.66,1.33,...,3.66}
        \fill ($(\x,0) + {\x}*(0,-0.375)$) circle (2pt);

      \fill ($(1.33,0) + 1.33*(0,-0.375)$) circle (3pt) node [above right] {$W_2$};

    \end{scope}
  \end{tikzpicture}
  \end{center}

  \caption{Illustration of the proof of Lemma~\ref{lemma:bad-event};
  the situation in $G$ is on the left, and the situation in $G'$ is on the right.}
  \label{figure:bad-event}
\end{figure}
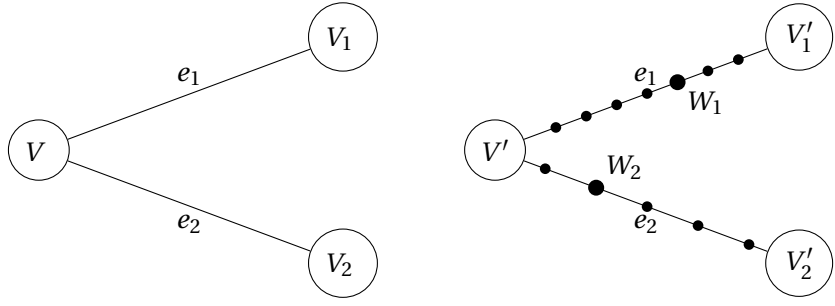
\begin{lemma}
  \label{lemma:bad-event}
  $\Pr[B_{e_1 e_2}] \leq \frac{4}{r}$.
\end{lemma}
\begin{proof}
  Let $\pi_1$ and $\pi_2$ be the permutations of rows chosen for $e_1$ and $e_2$, respectively,
  from the perspective of $V'$.
  That is, when walking along the subdivision path of $e_j$ in $G'$
  from $V'$ towards the opposite end-point, which we call $V_j$,
  then the first segment exchanges the symbols of the $\pi_j(1)$-th row.
  We will show that if $\pi_1(1)$, $\pi_1(2)$, $\pi_2(1)$ and $\pi_2(2)$ are all distinct,
  then the bad event $B_{e_1 e_2}$ does not occur.
  By the union bound, the probability of $B_{e_1 e_2}$ is therefore at most $4/r$.

  So suppose $\pi_1(1)$, $\pi_1(2)$, $\pi_2(1)$ and $\pi_2(2)$ are all distinct,
  and let $W_1 = \{ B_1 \} \in \cV'$ and $W_2 = \{ B_2 \} \in \cV'$ be two intermediate vertices
  on the subdivisions of $e_1$ and $e_2$ respectively,
  see Figure~\ref{figure:bad-event}.
  Let $W_1$ lie on the $i$-th segment of the subdivision path of $e_1$,
  that is, the $\pi_1(i)$-th row of $B_1$ is the only one that may differ from both $A \in V$ and $A_1 \in V_1$.
  Similarly, let $W_2$ lie on the $j$-th segment of the subdivision path of $e_2$.
  We claim that $B_1$ and $B_2$ differ in at least two rows.

  If both $i, j \leq 2$,
  then only the $\pi_1(1)$-th and $\pi_1(2)$-th row of $B_1$ may differ from $A$,
  and therefore $B_1|_{\pi_2(1)} = A$.
  On the other hand, the $\pi_2(1)$-th row of $B_2$ is different from $A$ by at least one element,
  and therefore $B_2|_{\pi_2(1)} \neq A$.
  Conversely, we have $B_1|_{\pi_1(1)} \neq A$ and $B_2|_{\pi_2(1)} = A$.
  So $B_1$ and $B_2$ differ in at least two rows.

  If both $i, j \geq 3$,
  then $\pi_1(i)$ is different from at least one of $\pi_2(1), \pi_2(2)$.
  Suppose that $\pi_1(i) \neq \pi_2(1)$.
  Then $B_1|_{\pi_2(1)} \in \{ A, A_1 \}$, whereas $B_2|_{\pi_2(1)} = A_2$.
  The other case is analogous,
  so that $B_1$ and $B_2$ differ in at least one of the $\pi_2(1)$-th and $\pi_2(2)$-th row.
  Similarly, they differ in at least one of the $\pi_1(1)$-th and the $\pi_1(2)$-th row,
  and so they differ in at least two rows..

  If $i \geq 3$ and $j \leq 2$,
  then $B_1|_{\pi_1(1)} = B_1|_{\pi_1(2)} = A_1$,
  while $B_2|_{\pi_1(1)} = B_2|_{\pi_1(2)} = A$,
  so that $B_1$ and $B_2$ again differ in at least two rows.
  The last case is symmetric.

  Since $B_1$ and $B_2$ differ in at least two rows in every case,
  we have $|B_1 \cap B_2| \leq rd - 2$.
  This applies to all pairs of vertices on the subdivision paths of $e_1$ and $e_2$ simultaneously,
  so the bad event $B_{e_1 e_2}$ does not occur
  if $\pi_1(1)$, $\pi_1(2)$, $\pi_2(1)$ and $\pi_2(2)$ are all distinct.
\end{proof}
We will now use Lovász' Local Lemma to prove that
with positive probability, none of the bad events happen. Recall~\cite{MR1885388}:
\begin{theorem}[Lovász' Local Lemma]
  Let $\cB = \{ B_i : i \in I \}$ be a finite set of events in a probability space
  such that $\Pr[B_i] \leq p$ for all $i \in I$.
  Furthermore, suppose that for every $i \in I$ there exists a neighborhood
  $N_i \subset \cB$, $|N_i| \leq k$,
  such that $B_i$ is independent from the set $B \setminus (N_i \cup \{ B_i \})$.
  If $(k+1)pe < 1$, then $\Pr[\mbox{none of the $B_i$ occur}] > 0$.
\end{theorem}
In our case, we can define the neighborhood of an event $B_{e_1 e_2}$
to be the set of all other bad events of the form $B_{e_1 f}$ and $B_{e_2 f}$.
Their number is bounded from above by $2\cdot (2\Delta - 3) = 4\Delta - 6$.
So by combining Lovász' Local Lemma and Lemma~\ref{lemma:bad-event},
it is sufficient to have
\[ (4\Delta - 5) \cdot \frac{4}{r} \cdot e < 1 \]
So it is sufficient to set $r = \lceil 16e\Delta \rceil$.
This finishes the proof of Lemma~\ref{lemma:adjacency},
and we can summarize the result of this section as follows:
\begin{theorem}
  \label{thm:construction}
  Let $G = (\cV, E)$ be a $d$-dimensional subset partition graph on the symbol set $S$
  that satisfies the singleton property.
  Let $r \geq \lceil 16e\Delta \rceil$, where $\Delta$ is the maximum degree of $G$.
  There exists an $rd$-dimensional SPG $G' = (\cV', E')$
  on the symbol set $S' = S \times [r]$ such that:
  \begin{enumerate}
    \item As an undirected graph, $G'$ is isomorphic to a subdivision of $G$.
    \item Let $V \in \cV$, $V = \{ A \}$, and let $V' \in \cV'$ the corresponding
      vertex under the mentioned isomorphism. Then $V' = \{ A \times [r] \}$.
    \item $G'$ satisfies strong adjacency, end-point count, and the singleton property.
  \end{enumerate}
\end{theorem}

\section{A spindle construction}

Let $d \geq 1$ and $n = 2d$.
Using $S = [d] \times \{ 1, 2 \}$ as symbol set,
one easily constructs an abstract spindle of exponential length.
Simply arrange all $d$-subsets of $S$ into an arbitrary order,
except that $A_1 = S \times \{1\}$ should come first,
and $A_2 = S \times \{2\}$ should come last.
We then construct a path $G$ whose vertices contain the $d$-subsets of $S$ in that order.
Then $G$ is a subset partition graph with the singleton property,
and moreover it is a spindle with apices $A_1$ and $A_2$
with length exponential in $d$.
By Theorem~\ref{thm:construction},
there exists a $32ed$-dimensional spindle $G'$ on $64ed$ symbols
of at least the same length
(in fact, its length will be greater by at least a factor of $32e$).
Moreover, $G'$ satisfies the strong adjacency and end-point count properties.

\section{Discussion}

There is some room for improving the constants in the bound on $r$ in Theorem~\ref{thm:construction},
but we did not explore this for sake of clarity.
The special case when the underlying graph is a path can be dealt with more easily
and with a better parameter $r$,
but we feel that further development of Theorem~\ref{thm:construction}
should go hand in hand with an investigation of which graph structures
are of interest in the study of subset partition graphs.

We consider the construction given in this paper
to be new evidence that the adjacency and end-point count properties are rather weak by themselves.
However, there remains the important question of
how they interact with the dimension reduction property.
Note that the construction of Theorem~\ref{thm:construction} preserves dimension reduction
only when the underlying structure of $G$ is a path.
If $G$ contains cycles, then dimension reduction is not preserved in general,
and fixing this problem would be a great step forward.

Right now, there is a gap between the best known constructions.
For SPGs that satisfy dimension reduction,
Eisenbrand et al.~\cite{EHRR10} give a construction with almost quadratic diameter.
However, their construction is far from satisfying the strong adjacency or end-point constructions.
On the other hand, the best known construction for SPGs
that satisfy both dimension reduction and strong adjacency
is via Santos'~\cite{Santos10} construction of polytopes with high (but still linear) diameter.
Closing this gap, perhaps by an improvement of the method used in this paper,
is an important open problem.

\small
\bibliographystyle{alpha}
\bibliography{books,mybib,papers,my_publications,local}

\end{document}